\documentclass{article}

\usepackage{amsmath}
\usepackage{amssymb}

\newtheorem{thm}{Theorem}
\newtheorem{conj}{Conjecture}

\newtheorem{cor}{Corollary}
\newtheorem{ques}{Question}

\title{An almost all result on $q_1 q_2 \equiv c \pmod q$}
\author{Tsz Ho Chan}

\begin{document}
\maketitle
\begin{abstract}
In this paper we consider the congruence equation $q_1 q_2 \equiv c
\pmod q$ with $a < q_1 \leq a + q^{1/2+\epsilon}$ and $b < q_2 \leq
b + q^{1/2+\epsilon}$ and show that it has solution for almost all
$a$ and $b$. Then we apply it to a question of Fujii and Kitaoka as
well as generalize it to more variables. At the end, we will present
a new way to attack the above congruence equation question through
higher moments.
\end{abstract}

\section{Introduction and main results}

A famous congruence equation question is the following:
\begin{ques} \label{q1}
Given any $\epsilon > 0$. Is it true that, for any modulus $q \geq
1$ and any integer $c$ with $(c, q) = 1$, the congruence equation
$q_1 q_2 \equiv c \pmod q$ is solvable for some $1 \leq q_1, q_2
\ll_\epsilon q^{1/2 + \epsilon}$?
\end{ques}

Davenport [\ref{D}] used Kloosterman sum estimates to show that the
above question is true for all $\epsilon > 1/3$. Using Weil's bound
on Kloosterman sums (see equation (\ref{kloos})), Davenport's
argument implies the truth of Question \ref{q1} for all $\epsilon >
1/4$. Recently in [\ref{S1}], Shparlinski got the same result with
the further restriction that $q_1$, $q_2$ are relatively prime to
one another. When $q$ is a prime number, Garaev [\ref{G1}] obtained
a slight improvement that Question \ref{q1} is true for all
$\epsilon \geq 1/4$.

Question \ref{q1} seems to be hard. How about proving it for almost
all $c$? Recently Garaev and Karatsuba [\ref{GKa}], and Shparlinski
[\ref{S}] proved that the above question is true for almost all $c$
with any $\epsilon > 0$ when $q$ is prime and $q$ in general
respectively. Their results are more general as one of the interval
can be replaced by a sufficiently large subset of the interval and
the other interval does not have to start from $1$. Furthermore when
$q$ is prime, Garaev and Garcia [\ref{GG}] showed the above almost
all result with $q_1$, $q_2$ in any intervals of length $q^{1/2 +
\epsilon}$ by considering solutions to $q_1 q_2 \equiv q_3 q_4 \pmod
q$. It used both character sum technique of [\ref{A}] and
exponential sum technique of [\ref{G1}].

Thus, in general, one does not have to restrict the ranges of $q_1$
and $q_2$ to start from $1$. In fact, the above question should be
true for $q_1$ and $q_2$ in any interval of length $O_\epsilon
(q^{1/2 + \epsilon})$. In this paper, we will prove that this is
indeed the case for almost all such pairs of intervals for $q_1$ and
$q_2$, namely
\begin{thm} \label{thm1}
For any modulus $q \geq 1$ and any integers $1 \leq L \leq q$ and
$(c,q) = 1$,
$$S := \sum_{a = 1}^{q} \sum_{b = 1}^{q} \Bigl( \mathop{\sum_{q_1
\in (a, a+L]} \mathop{{\sum\nolimits}'}_{q_2 \in (b, b+L]}}_{q_1 q_2
\equiv c \pmod q} 1 - \frac{1}{q} \sum_{q_1 \in (a, a+L]}
\mathop{{\sum\nolimits}'}_{q_2 \in (b, b+L]} 1 \Bigr)^2 \ll L^2 q
d(q)^3$$ where $\sum'$ means summing over those numbers that are
relatively prime to $q$.
\end{thm}

Let us interpret Theorem \ref{thm1}. Since
$$\mathop{{\sum\nolimits}'}_{q_2 \in (b, b+L]} 1 = \frac{\phi(q)}{q}
L + O(d(q)),$$ we have
\begin{align*}
\frac{1}{q} \sum_{q_1 \in (a, a+L]} \mathop{{\sum\nolimits}'}_{q_2
\in (b, b+L]} 1 =& \frac{\phi(q)}{q^2} L^2 + O \Bigl( \frac{d(q)
L}{q} + \frac{L}{q^{1/2 - \epsilon/2}} \Bigr) \\
=& \frac{\phi(q)}{q^2} L^2 + O \Bigl( \frac{L}{q^{1/2 - \epsilon/2}}
\Bigr)
\end{align*}
by $d(q) \ll q^{\epsilon / 2}$. Note that the error term is smaller
than the main term when $q^{1/2 + \epsilon} \ll L$. Thus, if we let
$N$ be the number of pairs of $(a,b)$ such that $q_1 q_2 \equiv c
\pmod q$ has no solution with $q_1 \in (a, a+L]$ and $q_2 \in (b,
b+L]$, then
\[
N \Bigl(\frac{\phi(q)}{q^2} L^2\Bigr)^2 \ll L^2 q d(q)^3.
\]
This implies $N \ll q^{3 + \epsilon} / L^2 \leq q^{2 - \epsilon}$
using $\phi(q) \gg q / \log \log q$ and $L \gg q^{1/2 + \epsilon}$.
Consequently, with $L = C_\epsilon q^{1/2 + \epsilon}$ where
$C_\epsilon > 0$ is large enough, we have
\begin{cor} \label{cor0}
Given any modulus $q \geq 1$ and any integer $c$ with $(c, q) = 1$.
For any $\epsilon > 0$, there exists a constant $C_\epsilon > 0$
such that
\[
q_1 q_2 \equiv c \pmod q
\]
is solvable in $q_1 \in (a, a + C_\epsilon q^{1/2 + \epsilon}]$,
$q_2 \in (b, b + C_\epsilon q^{1/2 + \epsilon}]$ for almost all
pairs of $a$ and $b$.
\end{cor}

\bigskip

Interestingly, Theorem \ref{thm1} can recover the currently best
result to Question \ref{q1}:
\begin{cor} \label{cor1}
Given any modulus $q \geq 1$ and any integer $c$ with $(c, q) = 1$.
For any $\epsilon > 0$ and integers $a$ and $b$, there exists a
constant $C_\epsilon > 0$ such that
\[
q_1 q_2 \equiv c \pmod q
\]
is solvable for some $q_1 \in (a, a + C_\epsilon q^{3/4 +
\epsilon}]$, $q_2 \in (b, b + C_\epsilon q^{3/4 + \epsilon}]$.
\end{cor}
This opens up a new way to attack Question \ref{q1} through looking
at higher moment analogue of Theorems \ref{thm1} or \ref{thm3}. We
shall discuss this in the last section.

In a similar spirit, Fujii and Kitaoka [\ref{FK}] studied the
\begin{ques} \label{q2}
For any lattice point $(x,y) \in \mathbb{Z}^2$, let $C_{(x,y)}(r)$
denote the compact disc with center $(x,y)$ and radius $r$. Let $q$
be a large number. Find the infimum $r(q)$ of all real numbers $r$
such that the square $[0,q] \times [0,q]$ is covered by
$$ \mathop{\bigcup_{
x,y = 1}^{q - 1}}_{xy \equiv 1 \bmod q} C_{(x,y)}(r).
$$
\end{ques}

When $q$ is a prime number, they proved that $r(q) \ll q^{3/4} \log
q$ and conjectured that $r(q) \ll_\epsilon q^{1/2 + \epsilon}$ for
every $\epsilon > 0$. Garaev [\ref{G}] mentioned that the argument
of [\ref{GK}] gives $r(q) \ll q^{3/4}$ for prime $q$. Using Theorem
\ref{thm1}, we can answer Question \ref{q2} in an almost all sense:
\begin{cor} \label{cor2}
With the notations in Question \ref{q2}, let $\tilde{r}(q)$ be the
infimum of all real numbers $r$ such that
$$ \mathop{\bigcup_{
x,y = 1}^{q - 1}}_{xy \equiv 1 \bmod q} C_{(x,y)}(r)
$$
cover an area of $(1 + o(1)) q^2$ in the square $[0,q] \times
[0,q]$. Then $\tilde{r}(q) \ll_\epsilon q^{1/2 + \epsilon}$ for any
$\epsilon > 0$.
\end{cor}

By modifying the proof of Theorem \ref{thm1} slightly, one can get
\begin{thm} \label{thm2}
For any modulus $q \geq 1$ and any integers $1 \leq L_1, L_2 \leq q$
and $(c,q) = 1$,
$$\sum_{a = 1}^{q} \sum_{b = 1}^{q} \Bigl( \mathop{\sum_{q_1
\in (a, a+L_1]} \mathop{{\sum\nolimits}'}_{q_2 \in (b, b+L_2]}}_{q_1
q_2 \equiv c \pmod q} 1 - \frac{1}{q} \sum_{q_1 \in (a, a+L_1]}
\mathop{{\sum\nolimits}'}_{q_2 \in (b, b+L_2]} 1 \Bigr)^2 \ll L_1
L_2 q d(q)^3$$ where $\sum'$ means summing over those numbers that
are relatively prime to $q$.
\end{thm}
Then one can discuss the above results for rectangles and eclipses
instead of squares and circles. We leave these for the readers to
explore.

\bigskip

More generally, one can consider the more variable version:
\begin{equation} \label{tcong}
q_1 q_2 ... q_t \equiv c \pmod q.
\end{equation}
Using character sum method, Shparlinski and Winterhof [\ref{SW}]
recently proved that for any $\epsilon > 0$ and $(c,q) = 1$,
\[
q_1 q_2 q_3 \equiv c \pmod q \hbox{ is solvable for some } 1 \leq
q_1, q_2, q_3 \leq q^{2/3 + \epsilon}
\]
and, for $t \geq 4$,
\[
q_1 q_2 ... q_t \equiv c \pmod q \hbox{ is solvable for some } 1
\leq q_1, q_2, ..., q_t \leq q^{1/3 + 1/ (t+2) + \epsilon}.
\]

We shall imitate Theorem \ref{thm1} and prove
\begin{thm} \label{thm3}
For any modulus $q \geq 1$ and any integers $1 \leq L \leq q$ and
$(c,q) = 1$,
\begin{align*}
S := & \sum_{a_1 = 1}^{q} \cdots \sum_{a_t = 1}^{q} \Big|
\mathop{\mathop{{\sum\nolimits}'}_{q_1 \in (a_1, a_1 + L]} \cdots
\mathop{{\sum\nolimits}'}_{q_t \in (a_t, a_t + L]}}_{q_1 \cdots q_t
\equiv c (\bmod q)} 1 - (\frac{L}{q})^t \phi(q)^{t-1} \Big|^2 \\
\leq & C_q^2 t^{2 \omega(q)} q^{t-1} L^t d(q)^t \Bigl[1 + \frac{t
(L+1)^{t-1}}{q} \Bigr]
\end{align*}
where $C_q = 1$ if $q$ is odd and $C_q = 2^{(t+1)/2}$ if $q$ is
even.
\end{thm}
\begin{cor} \label{cor3}
Given any modulus $q \geq 1$ and any integer $c$ with $(c, q) = 1$.
For any $\epsilon > 0$, there exists a constant $C_\epsilon > 0$
such that
\[
q_1 q_2 ... q_t \equiv c \pmod q
\]
is solvable in $q_1 \in (a_1, a_1 + C_\epsilon q^{1/t + \epsilon}]$,
$q_2 \in (a_2, a_2 + C_\epsilon q^{1/t + \epsilon}]$, ..., $q_t \in
(a_t, a_t + C_\epsilon q^{1/t + \epsilon}]$ for almost all
$t$-tuples $a_1$, ..., $a_t$.
\end{cor}
The exponent $1/t + \epsilon$ is best possible. One may then imitate
Corollary \ref{cor1} and get a non-almost all result for
(\ref{tcong}). By averaging over $c$, one can show that $S \gg
q^{t-1} L^t$ for some $c$. Thus, even with the best possible upper
bound $O(q^{t-1} L^t)$ for Theorem \ref{thm3}, one can only prove
that (\ref{tcong}) has solution for $q_1, ... , q_t$ in intervals of
length $q^{(t+1)/(2t)}$. These are no better than Shparlinski and
Winterhof's results. So passing from our almost all result to
non-almost all result is not a good approach unless one considers
higher moments or can somehow generate more tuples of intervals
without a solution out of a single one.

\bigskip

In summary, the method to study the above questions falls into two
categories. One uses exponential sums, particularly Kloosterman and
hyper-Kloosterman sums. The other one is character sum techniques
including Polya-Vinogradov and Burgess bounds as well as fourth
moment estimates on character sums (see [\ref{FI}] and [\ref{A}]).
It seems that character sum does better when there are more
variables. However, for our almost all results, we shall use
Kloosterman and hyper-Kloosterman sums.

\bigskip

The paper is organized as follows. First we will prove Theorem
\ref{thm1}. The reason we do this first is that it is how this
research began and it illustrates the essence of techniques used.
Then we will prove Corollaries \ref{cor1} and \ref{cor2}. After
these, we will prove the general case, Theorem \ref{thm3} and
Corollary \ref{cor3}, more neatly using the language of finite
Fourier series. Finally we will discuss higher moment attack of
Question \ref{q1}.

\bigskip

{\bf Notations} Throughout the paper, $\epsilon$ denotes a small
positive number. $f(x) \ll g(x)$ means that $|f(x)| \leq C g(x)$ for
some constant $C > 0$ and $f(x) \ll_\lambda g(x)$ means that the
implicit constant $C = C_\lambda$ may depend on the parameter
$\lambda$. Also $\phi(n)$ is Euler's phi function, $d(n)$ is the
number of divisors of $n$ and $\omega(n)$ is the number of distinct
prime divisors of $n$.
\section{Theorem \ref{thm1}}
For $(q_2, q) =1$, the congruence equation $q_1 q_2 \equiv c \pmod
q$ is equivalent to $q_1 \equiv c \overline{q_2} \pmod q$ where
$\overline{n}$ denotes the multiplicative inverse of $n \pmod q$. By
the orthogonal property of $e(u) = e^{2 \pi i u}$,
\[
\frac{1}{q} \sum_{a = 1}^{q} e\Bigl( \frac{k a}{q} \Bigr) = \left\{
\begin{tabular}{ll} $1$, & if $k \equiv 0 \bmod q$, \\
$0$, & otherwise,
\end{tabular} \right.
\]
we have
$$S = \sum_{a = 1}^{q} \sum_{b = 1}^{q} \Bigl( \frac{1}{q} \sum_{k=1}^{q-1}
\sum_{q_1 \in (a, a+L]} \mathop{{\sum\nolimits}'}_{q_2 \in (b, b+L]}
e \Bigl(\frac{k (q_1 - c \overline{q_2})}{q} \Bigr) \Bigr)^2.$$
Expanding things out, we have
\begin{align*}
S =& \sum_{a = 1}^{q} \sum_{b = 1}^{q} \frac{1}{q^2}
\sum_{k=1}^{q-1} \sum_{l=1}^{q-1} \sum_{q_1 \in (a, a+L]}
\mathop{{\sum\nolimits}'}_{q_2 \in (b, b+L]} e \Bigl(\frac{k (q_1 -
c \overline{q_2})}{q} \Bigr) \\
& \times \sum_{q_3 \in (a, a+L]} \mathop{{\sum\nolimits}'}_{q_4 \in
(b, b+L]} e \Bigl(\frac{-l (q_3 - c \overline{q_4})}{q} \Bigr) \\
=& \sum_{a = 1}^{q} \sum_{b = 1}^{q} \frac{1}{q^2} \sum_{k=1}^{q-1}
\sum_{l=1}^{q-1} \sum_{q_1, q_3 \in (a, a+L]}
\mathop{{\sum\nolimits}'}_{q_2, q_4 \in (b, b+L]} e\Bigl(\frac{k
q_1}{q}\Bigr) e\Bigl(\frac{-l q_3}{q}\Bigr) e\Bigl(\frac{-k c
\overline{q_2}}{q}\Bigr) e\Bigl(\frac{l c \overline{q_4}}{q}\Bigr).
\end{align*}
By making a change of variable $q_1 = q_3 + s$ with $-L \leq s \leq
L$ and combining the sums over $q_1$, $q_3$ and $a$, we have
\begin{align*}
S =& \frac{1}{q^2} \sum_{k=1}^{q-1} \sum_{l=1}^{q-1} \sum_{q_3 =
1}^{q} e\Bigl(\frac{(k-l) q_3}{q}\Bigr) \sum_{s = -L}^{L} (L - |s|)
e\Bigl(\frac{k s}{q} \Bigr) \\
&\times \sum_{b = 1}^{q} \mathop{{\sum\nolimits}'}_{q_2, q_4 \in (b,
b+L]} e\Bigl(\frac{-k c \overline{q_2}}{q}\Bigr) e\Bigl(\frac{l c
\overline{q_4}}{q}\Bigr) \\
=& \frac{1}{q} \sum_{k=1}^{q-1} \sum_{s = -L}^{L} (L - |s|)
e\Bigl(\frac{k s}{q} \Bigr) \sum_{b=1}^{q}
\mathop{{\sum\nolimits}'}_{q_2 \in (b, b+L]} e\Bigl(\frac{-k c
\overline{q_2}}{q}\Bigr) \mathop{{\sum\nolimits}'}_{q_4 \in (b,
b+L]} e\Bigl(\frac{k c \overline{q_4}}{q}\Bigr) \\
=& \frac{1}{q} \sum_{k=1}^{q-1} \Bigl(\frac{ \sin
\frac{Lk\pi}{q}}{\sin \frac{k\pi}{q}} \Bigr)^2 \sum_{b=1}^{q}
\mathop{{\sum\nolimits}'}_{q_2 \in (b, b+L]} e\Bigl(\frac{-k c
\overline{q_2}}{q}\Bigr) \mathop{{\sum\nolimits}'}_{q_4 \in (b,
b+L]} e\Bigl(\frac{k c \overline{q_4}}{q}\Bigr)
\end{align*}
by Fej\'{e}r kernel formula $(\frac{ \sin \pi N x}{\sin \pi x})^2 =
\sum_{j = -N}^{N} (N - |j|) e(j x)$. The innermost sums are
incomplete Kloosterman sums. We can use standard technique to
complete the sums:
\begin{align*}
\mathop{{\sum\nolimits}'}_{q_2 \in (b, b+L]} e\Bigl(\frac{-k c
\overline{q_2}}{q} \Bigr) =& \mathop{{\sum\nolimits}'}_{q_2 = 1}^{q}
e\Bigl(\frac{-k c \overline{q_2}}{q}\Bigr) \sum_{m \in (b, b+L]}
\frac{1}{q} \sum_{l = 1}^{q} e \Bigl(\frac{l (m - q_2)}{q} \Bigr) \\
=& \frac{1}{q} \sum_{l = 1}^{q} \sum_{m \in (b, b+L]} e\Bigl(\frac{l
m}{q} \Bigr) \mathop{{\sum\nolimits}'}_{q_2 = 1}^{q} e\Bigl(\frac{-l
q_2 -k c \overline{q_2}}{q}\Bigr) \\
=& \frac{1}{q} \sum_{l = 1}^{q} \Bigl[ \sum_{m \in (b, b+L]}
e\Bigl(\frac{l m}{q} \Bigr) \Bigr] S(-l, -k c; q)
\end{align*}
where $S(a,b;q) := \sum'_{1 \leq n \leq q} e(\frac{a n + b
\overline{n}}{q})$ is the Kloosterman sum. Therefore,
\begin{align*}
S = & \frac{1}{q^3} \sum_{k=1}^{q-1} \Bigl(\frac{ \sin
\frac{Lk\pi}{q}}{\sin \frac{k\pi}{q}} \Bigr)^2 \sum_{l_1 = 1}^{q}
\sum_{l_2 = 1}^{q} \sum_{b = 1}^{q} \Bigl[ \sum_{m_1 \in (b, b+L]} e
\Bigl(\frac{l_1 m_1}{q} \Bigr) \Bigr] \\
& \times \Bigl[ \sum_{m_2 \in (b, b+L]} e \Bigl(\frac{-l_2 m_2}{q}
\Bigr) \Bigr] S(-l_1, -k c; q) \overline{S(-l_2, -k c; q)}.
\end{align*}
By making a change of variable $m_1 = m_2 + d$ with $-L \leq d \leq
L$, the sums over $m_1$, $m_2$ and $b$ combine to give
\begin{align*}
S =& \frac{1}{q^3} \sum_{k=1}^{q-1} \Bigl(\frac{\sin
\frac{Lk\pi}{q}}{\sin \frac{k\pi}{q}} \Bigr)^2 \sum_{l_1 = 1}^{q}
\sum_{l_2 = 1}^{q} \sum_{d = -L}^{L} (L - |d|) e\Bigl(\frac{l_1
d}{q} \Bigr) \sum_{m_2 = 1}^{q} e\Bigl(\frac{(l_1 - l_2) m_2}{q}
\Bigr) \\
& \times S(-l_1, -k c; q) \overline{S(-l_2, -k c; q)} \\
=& \frac{1}{q^2} \sum_{k=1}^{q-1} \Bigl(\frac{\sin
\frac{Lk\pi}{q}}{\sin \frac{k\pi}{q}} \Bigr)^2 \sum_{l_1 = 1}^{q}
\sum_{d = -L}^{L} (L - |d|) e\Bigl(\frac{l_1 d}{q} \Bigr) S(-l_1,
-k c; q) \overline{S(-l_1, -k c; q)} \\
=& \frac{1}{q^2} \sum_{k=1}^{q-1} \Bigl(\frac{\sin
\frac{Lk\pi}{q}}{\sin \frac{k\pi}{q}} \Bigr)^2 \sum_{l_1 = 1}^{q -
1} |S(-l_1, -k c; q)|^2 \Bigl(\frac{\sin \frac{L l_1 \pi}{q}}{\sin
\frac{l_1 \pi}{q}} \Bigr)^2 \\
&+ \frac{L^2}{q^2} \sum_{k=1}^{q-1} \Bigl(\frac{\sin
\frac{Lk\pi}{q}}{\sin \frac{k\pi}{q}} \Bigr)^2 |S(0, -k c; q)|^2.
\end{align*}
Now recall Weil's bound on Kloosterman sums (see [\ref{IK},
Corollary 11.12] for example)
\begin{equation} \label{kloos}
S(a,b;q) \ll (a,b,q)^{1/2} q^{1/2} d(q)
\end{equation}
and
\begin{align*}
T(N; q) :=& \sum_{k = 1}^{q-1} \Bigl(\frac{\sin \pi N k /q}{\sin \pi
k /q} \Bigr)^2 \\
=& \sum_{k = 1}^{q-1} \sum_{d = -N}^{N} (N - |d|) e
\Bigl( \frac{d k}{q} \Bigr) = \sum_{d = -N}^{N} (N - |d|) \sum_{k =
1}^{q-1} e \Bigl( \frac{d k}{q} \Bigr) \\
=& \sum_{d = -N}^{N} (N - |d|) \sum_{k = 1}^{q} e \Bigl( \frac{d
k}{q} \Bigr) - N^2 = q N - N^2
\end{align*}
for $0 \leq N \leq q$. As $T(N \pm q; q) = T(N; q)$, we have $|T(N;
q)| \leq q N$ for all $N \geq 0$. Also, by grouping the sum
according to the greatest common divisor of $k$ and $q$,
\begin{equation} \label{kk}
\sum_{k = 1}^{q-1} \Bigl(\frac{\sin \pi L k /q}{\sin \pi k /q}
\Bigr)^2 (k,q) = \sum_{d | q} d \mathop{\sum_{k' = 1}^{q/d -
1}}_{(k', q/d) = 1} \Bigl(\frac{\sin \frac{L k' \pi}{q/d}}{\sin
\frac{k' \pi}{q/d}} \Bigr)^2 \leq \sum_{d | q} d \frac{q}{d} L = q L
d(q).
\end{equation}
Since $(c,q) = 1$, using (\ref{kloos}) and (\ref{kk}), we have
\begin{align*}
S \ll& \frac{d(q)^2}{q} \sum_{k=1}^{q-1} \Bigl(\frac{\sin \frac{Lk
\pi}{q}}{\sin \frac{k\pi}{q}} \Bigr)^2 \sum_{l_1 = 1}^{q - 1} (l_1,
q) \Bigl(\frac{\sin \frac{L l_1 \pi}{q}}{\sin \frac{l_1 \pi}{q}}
\Bigr)^2 + \frac{L^2}{q} d(q)^2 q L d(q) \\
\leq & L d(q)^3 \sum_{k=1}^{q-1} \Bigl(\frac{\sin \frac{Lk
\pi}{q}}{\sin \frac{k\pi}{q}} \Bigr)^2 + L^3 d(q)^3 \ll L^2 q d(q)^3
\end{align*}
as $L \leq q$. This proves Theorem \ref{thm1}.
\section{Corollaries \ref{cor1} and \ref{cor2}}
\

Proof of Corollary \ref{cor1}: Let $L = [C_\epsilon q^{3/4 +
\epsilon}]$. Suppose there are some integers $a_0$ and $b_0$ such
that the congruence equation $q_1 q_2 \equiv c \pmod q$ has no
solution with $a_0 < q_1 \leq a_0 + 2 L$, $b_0 < q_2 \leq b_0 + 2
L$. Then by Theorem \ref{thm1}, we have
\[
L^2 \Bigl(0 - \frac{L^2}{q} \Bigr)^2 \ll L^2 q d(q)^3
\]
as the congruence equation $q_1 q_2 \equiv c \pmod q$ has no
solution with $a < q_1 \leq a + L$, $b < q_2 \leq b + L$ for all
$a_0 \leq a \leq a_0 + L$ and $b_0 \leq b \leq b_0 + L$. The above
inequality gives
\[
\frac{L^6}{q^2} \ll L^2 q d(q)^3 \; \hbox{ or } \; L^4 \ll q^3
d(q)^3.
\]
This leads to $[ C_\epsilon q^{3/4 + \epsilon} ] = L \ll q^{3/4}
d(q)^{3/4}$ which is impossible if $C_\epsilon$ is large enough
(using $d(q) \ll_\epsilon q^\epsilon$). Hence we have Corollary
\ref{cor1}.

\bigskip

Proof of Corollary \ref{cor2}: Set $L := [ C_\epsilon q^{1/2 +
\epsilon} ]$. For $1 \leq x, y \leq q - 1$ and $x y \equiv 1 \bmod
q$, define the square and circle centered at $(x,y)$ by
\[
S_{(x,y)}(L) := \{ (a,b) : x - L - 1 \leq a \leq x + L + 1, y - L -
1 \leq b \leq y + L + 1 \},
\]
and
\[
C_{(x,y)}(L) := \{ (a,b) : (a - x)^2 + (b - y)^2 \leq L^2 \}.
\]
If $\bigcup_{ x,y = 1, \; xy \equiv 1 \bmod q}^{q - 1} S_{(x,y)}(L)$
covers the square $[0,q] \times [0,q]$, then the circles
\\
$\bigcup_{ x,y = 1, \; xy \equiv 1 \bmod q}^{q - 1}
C_{(x,y)}(\sqrt{2}(L+1))$ would cover the square $[0,q] \times
[0,q]$ and we are done.

\bigskip

Consider $L < a , b < q - L$ (this is to avoid ``wrap" around
squares $\pmod q$ when applying Theorem \ref{thm1}). Suppose $(a,b)$
is not covered by $\bigcup_{ x,y = 1, \; xy \equiv 1 \bmod q}^{q -
1} S_{(x,y)}(L)$. Then the square $([a], [a]+L] \times ([b], [b]+L]$
does not contain any solution to $q_1 q_2 \equiv 1 \bmod q$ with
$q_1 \in ([a], [a]+L]$ and $q_2 \in ([b], [b]+L]$ for otherwise the
square $S_{(q_1, q_2)}(L)$ would contain $(a,b)$. Call $([a], [b])$
a ``bad" lattice point and let $N$ be the number of such ``bad"
lattice points. Theorem \ref{thm1} tells us that $N \frac{L^4}{q^2}
\ll L^2 q d(q)^3$. Hence $N \ll \frac{q^3 d(q)^3}{L^2} = o(q^2)$ for
$C_\epsilon$ large enough (using $d(q) \ll_\epsilon
q^{\epsilon/2}$). Since every $(a,b)$ with $L < a , b < q - L$ not
covered is associated to some ``bad" lattice points, the area not
covered by $\bigcup_{ x,y = 1, \; xy \equiv 1 \bmod q}^{q - 1}
S_{(x,y)}(L)$ and hence $\bigcup_{ x,y = 1, \; xy \equiv 1 \bmod
q}^{q - 1} C_{(x,y)}(\sqrt{2}(L+1))$ must be $o(q^2)$ and we have
Corollary \ref{cor2} since the area outside of $L < a , b < q - L$
is $O(q L) = o(q^2)$.
\section{Theorem \ref{thm3} and Corollary \ref{cor3}}

Proof of Theorem \ref{thm3}: Let
\[
\chi_a(m) = \left\{ \begin{array}{ll} 1 & \hbox{ if } a < m
\leq a+L (\bmod q), \\
0 & \hbox{ otherwise.}
\end{array} \right.
\]
and its Fourier coefficients
\[
\hat{\chi}_a(n) = \frac{1}{q} \sum_{l (\bmod q)} \chi_a(l) e\Bigl( -
\frac{l n}{q}\Bigr) = \frac{1}{q} \sum_{a < l \leq a + L} e\Bigl(-
\frac{l n}{q}\Bigr).
\]
Then $\chi_a(m) = \sum_{k (\bmod q)} \hat{\chi}_a(k) e(\frac{k
m}{q})$ as its finite Fourier series.

For $(c, q) = 1$, we have
\[
S = \sum_{a_1 = 1}^{q} \cdots \sum_{a_t = 1}^{q} \Big|
\mathop{\mathop{{\sum\nolimits}'}_{q_1 = 1}^{q} \cdots
\mathop{{\sum\nolimits}'}_{q_t = 1}^{q}}_{q_1 \cdots q_t \equiv c
(\bmod q)} \chi_{a_1}(q_1) \cdots \chi_{a_t}(q_t) - \Bigl(
\frac{L}{q} \Bigr)^t \phi(q)^{t-1} \Big|^2.
\]
Since $\hat{\chi}_{a}(q) = L / q$,
\[
\mathop{\mathop{{\sum\nolimits}'}_{q_1 = 1}^{q} \cdots
\mathop{{\sum\nolimits}'}_{q_t = 1}^{q}}_{q_1 \cdots q_t \equiv c
(\bmod q)} \hat{\chi}_{a_1}(q) \cdots \hat{\chi}_{a_t}(q) =
\Bigl(\frac{L}{q} \Bigr)^t \phi(q)^{t-1}.
\]
Thus by using the finite Fourier series of $\chi_a(m)$, we have
\[
S = \sum_{a_1 = 1}^{q} \cdots \sum_{a_t = 1}^{q} \Big|
\mathop{\mathop{{\sum\nolimits}'}_{q_1 = 1}^{q} \cdots
\mathop{{\sum\nolimits}'}_{q_t = 1}^{q}}_{q_1 \cdots q_t \equiv c
(\bmod q)} \mathop{{\sum\nolimits}^*}_{1 \leq k_1, ..., k_t \leq q}
\hat{\chi}_{a_1}(k_1) e\Bigl(\frac{k_1 q_1}{q} \Bigr) \cdots
\hat{\chi}_{a_t}(k_t) e\Bigl(\frac{k_t q_t}{q} \Bigr) \Big|^2
\]
where $\sum'$ means summing over those numbers that are relatively
prime to $q$ and the $*$ means that we sum over all possible $k$'s
except $k_1 = \cdots = k_t = q$. Expanding things out, we have
\begin{align*}
S = & \sum_{a_1, ..., a_t} \mathop{\mathop{{\sum\nolimits}'}_{q_1,
..., q_t}}_{q_1 \cdots q_t \equiv c (\bmod q)}
\mathop{\mathop{{\sum\nolimits}'}_{q_1', ..., q_t'}}_{q_1' \cdots
q_t' \equiv c (\bmod q)} \mathop{{\sum\nolimits}^*}_{k_1, ..., k_t}
\mathop{{\sum\nolimits}^*}_{k_1', ..., k_t'} \\
& \hat{\chi}_{a_1}(k_1) \overline{\hat{\chi}_{a_1}(k_1')}
e\Bigl(\frac{k_1 q_1}{q}\Bigr) e\Bigl(- \frac{k_1' q_1'}{q}\Bigr)
\cdots \hat{\chi}_{a_t}(k_t) \overline{\hat{\chi}_{a_t}(k_t')}
e\Bigl(\frac{k_t q_t}{q}\Bigr) e\Bigl(- \frac{k_t' q_t'}{q}\Bigr).
\end{align*}
Observe that
\begin{align*}
& \sum_{a_1 = 1}^{q} \hat{\chi}_{a_1}(k_1)
\overline{\hat{\chi}_{a_1}(k_1')} \\
=& \frac{1}{q^2} \sum_{a_1 = 1}^{q} \sum_{a_1 < l_1 \leq a_1 + L}
e\Bigl(-\frac{l_1 k_1}{q}\Bigr) \sum_{a_1 < l_1' \leq a_1 + L}
e\Bigl(\frac{l_1' k_1'}{q}\Bigr) \\
=& \frac{1}{q^2} \sum_{a_1 = 1}^{q} \sum_{a_1 < l_1 \leq a_1 + L}
e\Bigl(\frac{l_1(k_1' - k_1)}{q}\Bigr) \sum_{a_1 - l_1 < s_1 \leq
a_1 + L - l_1} e\Bigl(\frac{s_1 k_1'}{q}\Bigr) \\
=& \frac{1}{q^2} \sum_{l_1 = 1}^{q} e\Bigl(\frac{l_1(k_1' -
k_1)}{q}\Bigr) \sum_{-L \leq s_1 \leq L} (L - |s_1|) e \Bigl(
\frac{s_1 k_1'}{q} \Bigr) \\
=& \left\{ \begin{array}{ll} \frac{1}{q} \sum_{-L \leq s_1 \leq L}
(L - |s_1|) e(\frac{s_1 k_1'}{q}), & \hbox{ if } k_1 \equiv k_1'
(\bmod q), \\
0, & \hbox{ otherwise}
\end{array} \right.
\end{align*}
by substituting $l_1' = l_1 + s_1$ and moving the sum over $a_1$
inside. Thus
\begin{align*}
S =& \frac{1}{q^t} \mathop{\mathop{{\sum\nolimits}'}_{q_1, ...,
q_t}}_{q_1 \cdots q_t \equiv c (\bmod q)}
\mathop{\mathop{{\sum\nolimits}'}_{q_1', ..., q_t'}}_{q_1' \cdots
q_t' \equiv c (\bmod q)} \mathop{{\sum\nolimits}^*}_{k_1, ..., k_t}
e\Bigl(\frac{k_1(q_1 - q_1')}{q}\Bigr) \cdots e\Bigl(
\frac{k_t(q_t - q_t')}{q}\Bigr) \\
& \times \Bigl(\frac{\sin \frac{\pi L k_1}{q}}{\sin \frac{\pi
k_1}{q}} \Bigr)^2 \cdots \Bigl(\frac{\sin \frac{\pi L k_t}{q}}{\sin
\frac{\pi k_t}{q}} \Bigr)^2 \\
=& \frac{1}{q^t} \mathop{{\sum\nolimits}^*}_{k_1, ..., k_t}
\Bigl(\frac{\sin \frac{\pi L k_1}{q}}{\sin \frac{\pi k_1}{q}}
\Bigr)^2 \cdots \Bigl(\frac{\sin \frac{\pi L k_t}{q}}{\sin \frac{\pi
k_t}{q}} \Bigr)^2 \Big| \mathop{\mathop{{\sum\nolimits}'}_{q_1, ...,
q_t}}_{q_1 \cdots q_t \equiv c (\bmod q)} e \Bigl( \frac{k_1 q_1 +
... + k_t q_t}{q} \Bigr) \Big|^2
\end{align*}
by Fej\'{e}r kernel formula $(\frac{ \sin \pi N x}{\sin \pi x})^2 =
\sum_{j = -N}^{N} (N - |j|) e(j x)$. Here we use the convention that
$(\frac{\sin \pi L k / q}{\sin \pi k / q})^2 = L^2$ if $k = q$. The
sum over the $q$'s is a hyper-Kloosterman sum. Now recall
Weinstein's version [\ref{W}] of Deligne's result on
hyper-Kloosterman sums:
\[
\Big| \mathop{\mathop{{\sum\nolimits}'}_{q_1, ..., q_t}}_{q_1 \cdots
q_t \equiv c (\bmod q)} e \Bigl( \frac{k_1 q_1 + ... + k_t q_t}{q}
\Bigr) \Big| \leq C_q q^{\frac{t-1}{2}} t^{\omega(q)} (k_1, k_t,
q)^{\frac{1}{2}} \cdots (k_{t-1}, k_t, q)^{\frac{1}{2}}
\]
where $C_q = 1$ if $q$ is odd and $C_q = 2^{(t+1)/2}$ if $q$ is
even, and $(a, b, c)$ stands for the greatest common divisor of $a$,
$b$ and $c$.

Using the above bound, we have
\[
S \leq \frac{C_q^2 t^{2 \omega(q)}}{q}
\mathop{{\sum\nolimits}^*}_{k_1, ..., k_t} \Bigl(\frac{\sin
\frac{\pi L k_1}{q}}{\sin \frac{\pi k_1}{q}} \Bigr)^2 \cdots
\Bigl(\frac{\sin \frac{\pi L k_t}{q}}{\sin \frac{\pi k_t}{q}}
\Bigr)^2 (k_1, k_t, q) \cdots (k_{t-1}, k_t, q).
\]
We estimate the above sum according to whether $i$ of the $k$'s are
equal to $q$ with $0 \leq i < t$. If $i = 0$, then it is bounded by
\[
\Bigl[ \sum_{k = 1}^{q-1} \Bigl(\frac{\sin \frac{\pi L k}{q}}{\sin
\frac{\pi k}{q}} \Bigr)^2 (k, q) \Bigr]^t \leq q^t L^t d(q)^t
\]
by (\ref{kk}). If $i > 0$, there are two cases depending on $k_t =
q$ or $k_t \neq q$.

When $k_t = q$, there are $i - 1$ of the $k$'s that can be $q$. So
we have the bound
\[
\binom{t-1}{i-1} q^{i-1} L^{2i} \Bigl[ \sum_{k = 1}^{q-1}
\Bigl(\frac{\sin \frac{\pi L k}{q}}{\sin \frac{\pi k}{q}} \Bigr)^2
(k, q) \Bigr]^{t-i} \leq \binom{t-1}{i-1} q^{t-1} L^{t + i} d(q)^{t}
\]
where $q^{i-1} L^{2i}$ comes from the $i-1$ such $k$'s together with
$k_t$.

When $k_t \neq q$, we have the bound
\begin{align*}
& \binom{t-1}{i} L^{2i} \Bigl[ \sum_{k = 1}^{q-1} \Bigl(\frac{\sin
\frac{\pi L k}{q}}{\sin \frac{\pi k}{q}} \Bigr)^2 (k, q)
\Bigr]^{t-i-1} \Bigl[ \sum_{k_t = 1}^{q-1} \Bigl(\frac{\sin
\frac{\pi L k_t}{q}}{\sin \frac{\pi k_t}{q}} \Bigr)^2 (k_t, q)^i
\Bigr] \\
\leq&  \binom{t-1}{i} L^{2i} (q L d(q))^{t - i - 1} q^i L d(q) \leq
\binom{t-1}{i} q^{t-1} L^{t + i} d(q)^{t}
\end{align*}
where $L^{2i}$ comes from the $i$ such $k$'s and they also
contribute $(k_t, q)^i$ to the sum over $k_t$. Combining the above
upper bounds, we have
\[
S \leq \frac{C_q^2 t^{2 \omega(q)}}{q} q^t L^t d(q)^t \Bigl[1 +
\frac{\sum_{i = 1}^{t-1} \binom{t}{i} L^i}{q} \Bigr] \leq C_q^2 t^{2
\omega(q)} q^{t-1} L^t d(q)^t \Bigl[1 + \frac{t (L+1)^{t-1}}{q}
\Bigr]
\]
by $(L+1)^t - L^t \leq t (L+1)^{t-1}$. This proves Theorem
\ref{thm3}.

\bigskip
Proof of Corollary \ref{cor3}: Recall
\[
S = \sum_{a_1 = 1}^{q} \cdots \sum_{a_t = 1}^{q} \Big|
\mathop{\mathop{{\sum\nolimits}'}_{q_1 = 1}^{q} \cdots
\mathop{{\sum\nolimits}'}_{q_t = 1}^{q}}_{q_1 \cdots q_t \equiv c
(\bmod q)} \chi_{a_1}(q_1) \cdots \chi_{a_t}(q_t) -
\Bigl(\frac{L}{q}\Bigr)^t \phi(q)^{t-1} \Big|^2
\]
Let $N$ be the number of tuples $(a_1, ..., a_t)$ such that $q_1
\cdots q_t \equiv c (\bmod q)$ has no solution with $a_1 < q_1 \leq
a_1 + L$, ..., $a_t < q_t \leq a_t + L$. Then by Theorem \ref{thm3},
\[
N \frac{L^{2t}}{q^2} (\frac{\phi(q)}{q})^{2t-2} \leq C_q^2 t^{2
\omega(q)} q^{t-1} L^t d(q)^t \Bigl[1 + \frac{t (L+1)^{t-1}}{q}
\Bigr].
\]
Since we are going to pick $L \leq q^{1/(t-1)} - 1$, we have
\[
N \leq C_q^2 (1 + t) t^{2 \omega(q)} d(q)^t (\frac{q}{\phi(q)})^{2t
- 2} q^{t} \frac{q}{L^t}
\]
which is $o(q^t)$ if $L = q^{1/t + \epsilon}$. This gives Corollary
\ref{cor3}.
\section{Higher moment attack}
In general, one expects that the error in
\[
\mathop{\mathop{{\sum\nolimits}'}_{q_1 \in (a_1, a_1 + L]} \cdots
\mathop{{\sum\nolimits}'}_{q_t \in (a_t, a_t + L]}}_{q_1 \cdots q_t
\equiv c (\bmod q)} 1 - (\frac{L}{q})^t \phi(q)^{t-1}
\]
is about the square root of the main term when $L \gg q^{1/t +
\epsilon}$. Focusing on $t = 2$, we expect
\[
\mathop{\mathop{{\sum\nolimits}'}_{q_1 \in (a_1, a_1 + L]}
\mathop{{\sum\nolimits}'}_{q_2 \in (a_2, a_2 + L]}}_{q_1 q_2 \equiv
c (\bmod q)} 1 - (\frac{L}{q})^2 \phi(q) \ll_\epsilon
\frac{L}{q^{1/2 - \epsilon}}.
\]
Raising to the $k$-th power and summing over $a_1$, $a_2$, we arrive
at the following
\begin{conj} \label{c1}
For any $\epsilon > 0$ and positive integer $k$,
\[
\sum_{a_1 = 1}^{q} \sum_{a_2 = 1}^{q} \Big|
\mathop{\mathop{{\sum\nolimits}'}_{q_1 \in (a_1, a_1 + L]}
\mathop{{\sum\nolimits}'}_{q_2 \in (a_2, a_2 + L]}}_{q_1 q_2 \equiv
c (\bmod q)} 1 - (\frac{L}{q})^2 \phi(q) \Big|^k \ll_{\epsilon, k}
\frac{L^k}{q^{k/2 - 2 - \epsilon}}
\]
for $L \gg q^{1/2 + \epsilon}$.
\end{conj}
Theorems \ref{thm1} and \ref{thm3} show that Conjecture \ref{c1} is
true for $k = 2$. Now we imitate the proof of Corollary \ref{cor1}.
Suppose there are integers $a_0$ and $b_0$ such that the congruence
equation $q_1 q_2 \equiv c \pmod q$ has no solution with $a_0 < q_1
\leq a_0 + 2 L$, $b_0 < q_2 \leq b_0 + 2 L$. Then by Conjecture
\ref{c1}, we have
\[
L^2 \Big| 0 - \frac{L^2}{q} \Big|^k \ll_{\epsilon, k}
\frac{L^k}{q^{k/2 - 2 - \epsilon}}
\]
as the congruence equation $q_1 q_2 \equiv c \pmod q$ has no
solution with $a < q_1 \leq a + L$, $b < q_2 \leq b + L$ for all
$a_0 \leq a \leq a_0 + L$ and $b_0 \leq b \leq b_0 + L$. The above
inequality gives
\[
L \ll_{\epsilon, k} q^{1/2 + 1/(k + 2) + \epsilon}.
\]
Consequently, if $L \gg_{\epsilon, k} q^{1/2 + 1/(k + 2) +
\epsilon}$, then $q_1 q_2 \equiv c \pmod q$ always has a solution
with $a_1 < q_1 \leq a_1 + L$, $a_2 < q_2 \leq a_2 + L$ for any
$a_1$, $a_2$.

\bigskip

In particular, if Conjecture \ref{c1} is true for $k = 3$ or $k =
4$, then Question \ref{q1} is true for all $\epsilon > 1/5$ or
$\epsilon > 1/6$ respectively. These are better than the currently
best result. In fact, if Conjecture \ref{c1} is true for arbitrarily
large $k$, we would settle Question \ref{q1} for all $\epsilon > 0$.
So the next challenge is to prove Conjecture \ref{c1} say for $k =
4$ even with a slightly larger upper bound. This would be a major
breakthrough!

\bigskip

{\bf Acknowledgements} The author would like to thank Professors
Stephen Choi and Kai-Man Tsang for stimulating discussions leading
to this work during a visit at the University of Hong Kong in the
summer of 2007.

Department of Mathematical Sciences \\
University of Memphis \\
Memphis, TN 38152 \\
U.S.A. \\
tchan@memphis.edu
\end{document}